\documentclass[11pt,letterpaper,reqno]{amsart}
\usepackage{amsmath,amsthm,amsfonts,amssymb,mathtools,algpseudocode}
\usepackage{comment,bookmark,bm,color,float,caption,subcaption}

\hypersetup{pdfstartview={FitH}}

\newtheorem{Theorem}{{\bf Theorem}}[section]

\newtheorem{Algorithm}[Theorem]{{\bf Algorithm}}
\newtheorem{Proposition}[Theorem]{{\bf Proposition}}
\newtheorem{Definition}[Theorem]{{\bf Definition}}

\numberwithin{equation}{section}

\newcommand{\C}{\mathbb{C}}
\newcommand{\R}{\mathbb{R}}
\newcommand{\N}{\mathbb{N}}
\newcommand{\calO}{\mbox{\Large $\mathcal{O}$}}
\newcommand{\calo}{\mathcal{O}}
\newcommand{\calC}{\mathcal{C}}
\newcommand{\calP}{\mathcal{P}}
\newcommand{\calS}{\mathcal{S}}
\newcommand{\off}{\textup{off}}
\newcommand{\diag}{\text{diag}}
\newcommand{\ve}{\text{vec}}
\newcommand{\calR}{\mathcal{R}}
\newcommand{\calJ}{\mathcal{J}}
\newcommand{\calN}{\mathcal{N}}

\begin{document}

\title[Convergence of the Eberlein diagonalization method]{Convergence of the Eberlein diagonalization method under the generalized serial pivot strategies}
\author{Erna Begovi\'{c}~Kova\v{c}}\thanks{\textsc{Erna Begovi\'{c} Kova\v{c},
University of Zagreb Faculty of Chemical Engineering and Technology, Maruli\'{c}ev trg 19, 10000 Zagreb, Croatia},
\texttt{ebegovic@fkit.hr}}
\author{Ana Perkovi\'{c}}\thanks{\textsc{Ana Perkovi\'{c},
University of Zagreb Faculty of Chemical Engineering and Technology, Maruli\'{c}ev trg 19, 10000 Zagreb, Croatia},
\texttt{aperkov@fkit.hr}}

\thanks{This work has been supported in part by Croatian Science Foundation under the project UIP-2019-04-5200.}
\date{\today}

\renewcommand{\subjclassname}{\textup{2020} Mathematics Subject Classification}

\begin{abstract}
The Eberlein method is a Jacobi-type process for solving the eigenvalue problem of an arbitrary matrix.
In each iteration two transformations are applied on the underlying matrix, a plane rotation and a non-unitary elementary transformation.
The paper studies the method under the broad class of generalized serial pivot strategies.
We prove the global convergence of the Eberlein method under the generalized serial pivot strategies with permutations and present several numerical examples.
\end{abstract}

\maketitle

\section{Introduction}

The Jacobi diagonalization method is the method of choice for solving the eigenproblem for dense symmetric matrices. Compared to the other state-of-the-art diagonalization methods, the main advantage of the Jacobi method is its high relative accuracy~\cite{DeVe92,Matejas09}. The method has been modified to deal with different matrix structures~\cite{FMM01,HSS14,MMMM09,MMT03,Mehl04} and to tackle various problems of numerical linear algebra~\cite{BGFass97,DV08-1,DV08-2,Mehl08}. Its convergence has been extensively studied. (See, e.g.,\@ \cite{Hari15,Masc95}.) One of the generalizations of the Jacobi method is known as the Eberlein method.

The Eberlein method, originally proposed in 1962~\cite{Eber62}, is a Jacobi-type process for solving the eigenvalue problem of an arbitrary matrix. It is one of the first efficient norm-reducing methods of this type. The iterative process on a general matrix $A\in\C^{n\times n}$ takes the form
\begin{equation}\label{Eber_process}
A^{(k+1)}=T_k^{-1}A^{(k)}T_k, \quad k\geq0,
\end{equation}
where $A^{(0)}=A$, and
$$T_k=R_kS_k$$
are non-singular elementary matrices. In particular, matrices $R_k$ are plane rotations and $S_k$ are non-orthogonal elementary matrices. The transformations $R_k$ are chosen to annihilate the pivot element of the matrix $(A^{(k)}+(A^{(k)})^*)/2$, while the transformations $S_k$ reduce the Frobenius norm of $A^{(k)}$.

In~\cite{Ves76}, Veselić studied a slightly altered Eberlein algorithm where in the $k$th step only one transformation is applied, either $R_k$ or $S_k$, but not both at the same time. He proved the convergence of this modified method under the classical Jacobi pivot strategy. Later, Hari~\cite{Hari82} proved the global convergence of the original method under the column/row cyclic pivot strategy on real matrices. In~\cite{PupHari99} Hari and Pupovci proved the convergence of the Eberlein method on complex matrices with the pivot strategies that are weakly equivalent to the row cyclic strategy. In the same paper authors considered the parallel method and proved its convergence under the pivot strategies that are weakly equivalent to the modulus strategy.

In this paper we expand the global convergence result for the Eberlein method to a significantly broader class of the cyclic pivot strategies --- generalized serial strategies with permutations, studied in~\cite{BePhD,BH17,BH21}. We consider the method in the form that is given in~\cite{PupHari99}. Our new result is the global convergence of the Eberlein method under the generalized serial pivot strategies with permutations. It is given in Theorem~\ref{theorem:sp}. We show that for an arbitrary $n\times n$ starting matrix $A^{(0)}$, $k\geq0$, the sequence $A^{(k)}$ converges to a block diagonal normal matrix. At the same time, the sequence $(A^{(k)}+(A^{(k)})^*)/2$ converges to a diagonal matrix $\diag(\mu_1,\mu_2,\ldots,\mu_n)$, where $\{\mu_1,\mu_2,\ldots,\mu_n\}$ are the real parts of the eigenvalues of $A$. 
Moreover, we present several numerical examples and discuss the cases of the unique and the multiple eigenvalues.

The paper is organized as follows. In Section~\ref{sec:eberlein} we describe the Eberlein method, its complex and real variant, while in Section~\ref{sec:strategies} we characterize the set of the pivot strategies that we work with. The main part of the paper is contained in Section~\ref{sec:cvg} where we prove the convergence of the method under the strategies from Section~\ref{sec:strategies}. Finally, in Section~\ref{sec:num} we present the results of our numerical tests.

\section{The Eberlein method}\label{sec:eberlein}

As it was mentioned in the introduction, there are several variations of the Eberlein method. The method can be applied to complex matrices using transformations $T_k\in\C^{n\times n}$, $k\geq0$, or one can observe the real method with $T_k\in\R^{n\times n}$, $k\geq0$. In this paper we mostly focus on the complex method. We describe it in Subsection~\ref{sec:complex}. In Subsection~\ref{sec:real} we outline the real case. We use the notation $\imath=\sqrt{-1}$. For a complex number $x$, $\text{Re}(x)$ stands for the real part of $x$ and $\text{Im}(x)$ stands for its imaginary part.

\subsection{The complex case}\label{sec:complex}

The Eberlein method is an iterative Jacobi-type method used to find the eigenvalues and eigenvectors of an arbitrary matrix $A\in\C^{n\times n}$. One iteration step of the method is given by the relation~\eqref{Eber_process}.
In the $k$th iteration, transformation $T_k$ is an elementary matrix that differs from the identity only in one of its $2\times 2$ principal submatrices $\widehat{T}_k$ determined by the pivot pair $(p(k),q(k))$,
\begin{equation*}
\widehat{T}_k=\begin{bmatrix}
    t_{p(k)p(k)}^{(k)}&t_{p(k)q(k)}^{(k)}\\
    t_{q(k)p(k)}^{(k)}&t_{q(k)q(k)}^{(k)}
\end{bmatrix}.
\end{equation*}
Matrix $T_k$ is set to be the product of two nonsingular matrices, a plane rotation $R_k$ and a non-unitary elementary matrix $S_k$. That is, $T_k=R_kS_k$.
Denote the $k$th pivot pair by $(p,q)=(p(k),q(k))$. The pivot pair is the same for both $R_k$ and $S_k$, and consequently for $T_k$.
In addition to $(p,q)$, matrices $R_k$ and $S_k$, depend on transformation angles $\alpha_k, \varphi_k$, and $\beta_k, \psi_k$, respectively. The pivot submatrix $\widehat{T}_k$ is equal to $\widehat{T}_k=\widehat{R}_k\widehat{S}_k\in\C^{2\times2}$, where
\begin{equation}\label{eq:RS}
    \widehat{R}_k=\begin{bmatrix}
        \cos \varphi_{k} & -e^{\imath\alpha_k} \sin \varphi_k\\
        e^{-\imath\alpha_k}\sin \varphi_k & \cos \varphi_k
    \end{bmatrix}, \quad
    \widehat{S}_k=\begin{bmatrix}
        \cosh\psi_k &  -\imath e^{\imath\beta_k}\sinh \psi_k\\
         \imath e^{-\imath\beta_k}\sinh \psi_k & \cosh \psi_k
    \end{bmatrix}.
\end{equation}

The process~\eqref{Eber_process} can be written with an intermediate step
\begin{equation*}
    \widetilde{A}^{(k)}=R_k^* A^{(k)} R_k, \quad A^{(k+1)}=S_k^{-1} \widetilde{A}^{(k)} S_k,\quad k\geq 0.
\end{equation*}
Let
\begin{align}
B^{(k)} & =\frac{1}{2}(A^{(k)}+A^{(k)*}), \label{eq:B} \\
\widetilde{B}^{(k)} & =R_k^* B^{(k)} R_k. \nonumber
\end{align}
Next, let $C$ be an operator defined by
\begin{equation}\label{eq:C}
C(A)=AA^*-A^*A.
\end{equation}
We denote $C(A^{(k)})=(c^{(k)}_{ij})$ and $C(\widetilde{A}^{(k)})=(\tilde{c}^{(k)}_{ij})$.
Obviously, $C(A)=0$, if and only if $A$ is a normal matrix.

The rotation $R_k$ is chosen such that the element of $B^{(k)}$ on position $(p,q)$ is annihilated.
The real number $\alpha_k$ and the angle $\varphi_k$ in~\eqref{eq:RS} are calculated from the following expressions,
\begin{align}
    \alpha_k & =\arg(b^{(k)}_{pq}), \label{eq:alpha} \\
    \tan 2\varphi_k & =\frac{2|b^{(k)}_{pq}|}{b^{(k)}_{pp}-b^{(k)}_{qq}}, \quad |\varphi_k|\leq \frac{\pi}{4}. \label{eq:2phi}
\end{align}
These formulas are the same as for the complex Jacobi method on Hermitian matrices.
Then, in order to get $\sin\varphi$ and $\cos\varphi$, we use formulas
\begin{align}
\tan\varphi_k=\frac{2|b_{pq}|\text{sign}(b_{pp}-b_{qq})}{|b_{pp}-b_{qq}|+\sqrt{|b_{pp}-b_{qq}|^2+4|b_{pq}|^2}}, \nonumber \\
\cos\varphi_k=\frac{1}{\sqrt{1+\tan^2\varphi_k}}, \ \sin\varphi_k=\frac{\tan\varphi_k}{\sqrt{1+\tan^2\varphi_k}}. \label{eq:phi}
\end{align}

On the other hand, $S_k$ is chosen to reduce the Frobenius norm of $A^{(k)}$. Set
$$\Delta_k:=\|A^{(k)}\|_F^2 -\|A^{(k+1)}\|_F^2.$$
Eberlein (\cite{Eber62}) proved that
\begin{align*}
\Delta_k & = \|\widetilde{A}^{(k)}\|_F^2 - \|A^{(k+1)}\|_F^2 \\
& =g_{pq}^{(k)}(1-\cosh 2\psi_k) -h_{pq}^{(k)}\sinh2\psi_k+\frac{1}{2}(|\Tilde{\xi}_{pq}^{(k)}|^2
+|\Tilde{d}_{pq}^{(k)}|^2)(1-\cosh4\psi_k) \\
& \qquad +\text{Im}(\Tilde{\xi}_{pq}^{(k)}\Tilde{d}_{pq}^{(k)*})\sinh 4\psi_k,
\end{align*}
where
\begin{align*}
    g_{pq}^{(k)}&=\sum_{\substack{i=1 \\ i\neq p,q}}^n |\Tilde{a}_{ip}^{(k)}|^2+|\Tilde{a}_{pi}^{(k)}|^2+|\Tilde{a}_{iq}^{(k)}|^2+|\Tilde{a}_{qi}^{(k)}|^2,\\
    h_{pq}^{(k)}&=-\text{Re}(l_{pq}^{(k)})\sin \beta_k + \text{Im}(l_{pq}^{(k)})\cos \beta_k,\\
    l_{pq}^{(k)}&=2 \sum_{\substack{i=1 \\ i\neq p,q}}^n (\Tilde{a}_{pi}^{(k)}\Tilde{a}_{qi}^{(k)*}-\Tilde{a}_{ip}^{(k)*}\Tilde{a}_{iq}^{(k)}) \\
    \Tilde{d}_{pq}^{(k)}&=\Tilde{a}_{pp}^{(k)}-\Tilde{a}_{qq}^{(k)},\\
    \Tilde{\xi}_{pq}^{(k)}&=(\Tilde{a}_{pq}^{(k)}+\Tilde{a}_{qp}^{(k)})\cos \beta_k-i(\Tilde{a}_{pq}^{(k)}-\Tilde{a}_{qp}^{(k)})\sin \beta_k.
\end{align*}
It is shown in~\cite{Eber62} that the choice of $\beta_k$ and $\psi_k$ such that
\begin{align}
\tan\beta_k & = -\frac{\text{Re}(\Tilde{c}_{pq})}{\text{Im}(\Tilde{c}_{pq})}, \label{eq:beta} \\
\tanh\psi_k & =\frac{1}{2}\frac{2\text{Im}(\Tilde{\xi}_{pq}^{(k)}\Tilde{d}_{pq}^{(k)*})-h_{pq}^{(k)}} {g_{pq}^{(k)}+2(|\Tilde{\xi}_{pq}^{(k)}|^2+|\Tilde{d}_{pq}^{(k)}|^2)}, \nonumber \\
& \cosh\psi_k =\frac{1}{\sqrt{1-\tanh^2\psi_k}}, \ \sinh\psi_k=\frac{\tanh\psi_k}{\sqrt{1-\tanh^2\psi_k}}, \label{eq:psi}
\end{align}
implies
\begin{equation}\label{eq:delta}
\Delta_k \geq\frac{1}{3}\frac{|\Tilde{c}_{pq}^{(k)}|^2}{\|A^{(k)}\|_F^2} \geq\frac{1}{3}\frac{|\Tilde{c}_{pq}^{(k)}|^2}{\|A\|_F^2}, \quad k\geq1.
\end{equation}

We summarize this procedure in Algorithm~\ref{agm:eberlein}.

\begin{Algorithm}\label{agm:eberlein}
\hrule\vspace{1ex}
\emph{Eberlein method}
\vspace{0.5ex}\hrule
\begin{algorithmic}
\State \textbf{Input:} $A\in \C^{n\times n}$
\State \textbf{Output: } matrix $A^{(k)}$
\State $A^{(0)}=A$; $T^{(0)}=I_n$;
\State $k=0$
\Repeat
\State Choose pivot pair $(p,q)$ according to the pivot strategy.
\State Find $\alpha_k$ using~\eqref{eq:alpha}, and $\sin\varphi_k$, $\cos\varphi_k$ using~\eqref{eq:phi}.
\State $\widetilde{A}^{(k)}=R_k^*{A}^{(k)}R_k$
\State Find $\beta_k$ using~\eqref{eq:beta}, and $\sin\psi_k$, $\cos\psi_k$ using~\eqref{eq:psi}.
\State ${A}^{(k+1)}=S_k^{-1}\widetilde{A}^{(k)}S_k$
\State $k=k+1$
\Until{convergence}
\end{algorithmic}
\hrule
\end{Algorithm}

One should keep in mind that it is not needed to formulate matrices $\widetilde{A}^{(k)}$ explicitly since both transformations $R_k$ and $S_k$ effect only the elements from the $p$th and $q$th row and column of $A^{(k)}$.

\subsection{The real case}\label{sec:real}

Suppose that $A$ is a real matrix and we wish for the iterates $A^{(k)}$ to stay real during the process~\eqref{Eber_process}. In order to satisfy this request we modify the complex algorithm. Firstly, we can take $\alpha_k=\pi$ and $\beta_k=\pi/2$. This implies
\begin{equation*}
    \widehat{R}_k=\begin{bmatrix}
        \cos \varphi_{k} & \sin \varphi_k\\
       -\sin \varphi_k & \cos \varphi_k
    \end{bmatrix}, \quad
    \widehat{S}_k=\begin{bmatrix}
        \cosh \psi_k &  \sinh \psi_k\\
         \sinh \psi_k & \cosh \psi_k
    \end{bmatrix}.
\end{equation*}
It remains to calculate the angles $\varphi$ and $\psi$.

As in the complex case, $\varphi_k$ is selected to annihilate the pivot element of $B^{(k)}$ while $\psi_k$ is chosen to reduce $\|A^{(k)}\|_F$.
The angle $\varphi_k$ is calculated from the relation similar to~\eqref{eq:2phi},
\begin{equation*}\label{eq:phi_real}
    \tan 2\varphi_k  =\frac{2b^{(k)}_{pq}}{b^{(k)}_{qq}-b^{(k)}_{pp}}, \quad |\varphi_k|\leq \frac{\pi}{4}.
\end{equation*}
Considering that $\beta_k=\pi/2$ and that all the elements of $A^{(k)}$ are real, the formula for $\psi_k$ is transformed into
$$\tanh\psi_k=\frac{\tilde{c}_{pq}^{(k)}}{g_{pq}^{(k)}+2((\tilde{e}_{pq}^{(k)})^2+(\tilde{d}_{pq}^{(k)})^2)},$$
where
$$\tilde{e}_{pq}^{(k)}=\tilde{a}_{pq}-\tilde{a}_{qp},$$
while $g_{pq}^{(k)}$ and $\tilde{d}_{pq}^{(k)}$ are the same as in the complex case.

\section{Generalized serial pivot strategies}\label{sec:strategies}

In each iteration $k$ of the Algorithm~\ref{agm:eberlein}, pivot position is selected according to the pivot strategy. In this section we describe a large class of pivot strategies that we work with --- generalized serial pivot strategies with permutations defined in~\cite{BH17}.

For an $n\times n$ matrix, possible pivot pairs are those from the upper triangle of the matrix, $\calP_n\coloneqq \{(i,j):1\leq i<j\leq n\}$.
A pivot strategy is any function $I\colon \N_0 \to \calP_n$, $\N_0=\{0,1,2,\ldots\}$. We work with \emph{cyclic pivot strategies}. Thus, we take $I$ as a periodic function with period $T=N\equiv \frac{n(n-1)}{2}$ and image $\calP_n$.

Pivot strategies are often easier understood using pivot orderings.
A cyclic strategy $I$ defines a sequence $\calo_I$ which is an ordering of $\calP_n$,
$$\calo_I=I(0), I(1), \ldots, I(N-1) \in \calO(\calP_n),$$
where $\calO(\calP_n)$ stands for the set of all finite sequences of elements from $\calP_n$, provided that each pair from $\calP_n$ appears at least once in every sequence.
An \emph{admissible transposition} in a pivot sequence $\calo$ is any transposition of two adjacent pivot pairs,
$$ (i_r,j_r),(i_{r+1},j_{r+1}) \to (i_{r+1},j_{r+1}),(i_r,j_r),$$
assuming that the sets $\{i_r,j_r\}$ and $\{i_{r+1},j_{r+1}\}$ are disjoint. Moreover, we need several equivalence relations on pivot orderings. (See, e.g.,~\cite{BH17}.)

\begin{Definition}\label{def:equiv}
Two pivot sequences $\calo$ and $\calo'$, where $\calo=(i_0,j_0),(i_1,j_1),\ldots,(i_r,j_r)$, are said to be
\begin{enumerate}
        \item[(i)] \textup{equivalent} $(\calo \sim \calo')$ if one can be obtained from the other by a finite set of admissible transpositions
        \item[(ii)] \textup{shift-equivalent} $(\calo \stackrel{s}{\sim} \calo')$ if $\calo=\left[ \calo_1,\calo_2 \right]$ and $\calo'=\left[\calo_2,\calo_1\right]$, where $\left[~,~\right]$ denotes the concatenation; the length of $\calo_1$ is called the shift length
        \item[(iii)] \textup{weak equivalent} $(\calo \stackrel{w}{\sim} \calo')$ if there exist $\calo_i\in\calO(\calS)$, $0\leq i \leq t,$ such that every two adjacent terms in the sequence $\calo=\calo_0, \calo_1,\ldots,\calo_t=\calo'$ are equivalent or shift-equivalent
        \item[(iv)] \textup{permutation equivalent}  $(\calo \stackrel{p}{\sim} \calo'$ or $\calo'=\calo(\textup{q}))$ if there is a permutation $\textup{q}$ of the set $\calS$ such that $\calo'=(\textup{q}(i_0),\textup{q}(j_0)),(\textup{q}(i_1),\textup{q}(j_1)),\ldots,(\textup{q}(i_r),\textup{q}(j_r))$
        \item[(v)] \textup{reverse} $(\calo'=\calo^{\gets})$ if $\calo'=(i_r,j_r),\ldots,(i_1,j_1),(i_0,j_0)$.
\end{enumerate}
Two pivot strategies $I_\calo$ and $I_{\calo'}$ are equivalent (shift-equivalent, weak equivalent, permutation equivalent, reverse) if the same is true for their corresponding pivot orderings $\calo$ and $\calo'$.
\end{Definition}

It is easy to see that, if $\calo \stackrel{w}{\sim} \calo'$, then there is a finite sequence $\calo=\calo_0, \calo_1,\ldots,$ $\calo_t=\calo'$ such that
\begin{equation}\label{can_form}
    \calo\sim\calo_1\stackrel{s}{\sim}\calo_2\sim\calo_3\stackrel{s}{\sim}\calo_4\ldots\calo' \quad \text{or} \quad \calo\stackrel{s}{\sim}\calo_1\sim\calo_2\stackrel{s}{\sim}\calo_3\sim\calo_4\ldots\calo'.
\end{equation}
The chain from~\eqref{can_form} that is connecting $\calo$ and $\calo'$ is in the \textit{canonical form}.

The most intuitive cyclic strategies are the \emph{row-cyclic}, $I_{row}=I_{\calo_{row}}$, and the \emph{column-cyclic strategy}, $I_{col}=I_{\calo_{col}}$, collectively named \emph{serial pivot strategies}.
Cyclic strategies that are equivalent to the serial pivot strategies are called \emph{wavefront strategies}.

Eberlein~\cite{Eber62} used the strategy where the pivot pair $(p,q)$ is chosen such that $4|c_{pq}^{(k)}|^2+(c_{pp}^{(k)}-c_{qq}^{(k)})^2$ is greater or equal to the average of all possible results for
$4|c_{ij}^{(k)}|^2+(c_{ii}^{(k)}-c_{jj}^{(k)})^2$, $1\leq i<j\leq n$.
Veselić~\cite{Ves76} used the classical Jacobi pivot strategy which takes the pivot pair that is the largest in the absolute value. Employing both of those strategies slows the algorithm down for large matrices.
Later, Hari~\cite{Hari82} proved the convergence for the real method under the wavefront strategies. In~\cite{PupHari99} Pupovci and Hari provided the convergence proof for the complex method using the parallel modulus strategy and the strategies that are weakly equivalent to it.

Now, let us describe the generalized serial pivot strategies with permutations. More details on these strategies can be found in ~\cite{BH17} and~\cite{BH21}.
For $l_1<l_2$, denote by $\Pi^{(l_1,l_2)}$ the set of all permutations of the set $\{l_1,l_1+1,l_1+2,\ldots,l_2\}$.
Let
\begin{align*}
    \calC_c^{(n)}=\Big\{ \calo \in \calO(\calP_n) \mid \calo= &(1,2),(\tau_3(1),3),(\tau_3(2),3),\ldots,(\tau_n(1),n),\ldots\\
    &\ldots,(\tau_n(n-1),n), \quad \tau_j\in\Pi^{(1,j-1)}, 3\leq j\leq n \Big\}.
\end{align*}
The orderings from $\calC_c^{(n)}$ go through the matrix column by column, starting from the second one, just like in the standard column strategy $I_{col}$. However, in each column pivot elements are chosen in some arbitrary order. If $\calo\in\calC_c^{(n)}$, then $\calo$ is called a column-wise ordering with permutations.
Likewise, the set of row-wise orderings with permutations is defined as
\begin{align*}
    \calC_r^{(n)}=\Big\{ \calo \in \calO(\calP_n) \mid \calo= &(n-1,n),(n-2,\tau_{n-2}(n-1)),(n-2,\tau_{n-2}(n)),\ldots\\
    &\ldots,(1,\tau_1(2)),\ldots, (1,\tau_1(n)) \quad \tau_i\in\Pi^{(i+1,n)}, 1\leq i\leq n-2 \Big\}.
\end{align*}
By employing these two sets of orderings and their reverses we define the set of \emph{serial orderings with permutations},
$$\calC_{sp}^{(n)}=\calC_c^{(n)}\cup \overleftarrow{\calC}_c^{(n)} \cup \calC_r^{(n)}\cup \overleftarrow{\calC}_r^{(n)}.$$
We expand the set $\calC_{sp}^{(n)}$ using the equivalence relations from the Definition~\ref{def:equiv}. Let
\begin{equation*}
    \calC^{(n)}_{sg}=\Big\{ \calo \in\calO(\calP_n) \mid \calo \stackrel{w}{\sim} \calo' \stackrel{p}{\sim}\calo'' \text{ or } \calo \stackrel{p}{\sim} \calo' \stackrel{w}{\sim}\calo'',  \calo''\in\calO_{sp}^{(n)}\Big\},
\end{equation*}
where $\calo'\in\calO(\calP_n)$. Strategies defined by orderings from $\calC^{(n)}_{sg}$ are called \textit{generalized serial pivot strategies with permutations}.

\section{Convergence of the Eberlein method}\label{sec:cvg}

In this section we prove that the iterative process~\eqref{Eber_process} converges under any pivot ordering $\calo\in\calC^{(n)}_{sg}$.
First we list several auxiliary results from the literature and their direct implications. We use the notation introduced in Section~\ref{sec:eberlein}.

\begin{itemize}
\item[(i)] (Eberlein~\cite{Eber62}) For $\|A^{(k)}\|_F^2$ we have
\begin{equation}\label{eq:eberlein}
\Delta_k=\|A^{(k)}\|_F^2-\|A^{(k+1)}\|_F^2 = \|\widetilde{A}^{(k)}\|_F^2-\|A^{(k+1)}\|_F^2 \geq0.
\end{equation}
\item[(ii)] Since the sequence $(\|A^{(k)}\|_F^2, \, k\geq0)$ is nonincreasing and bounded, it is convergent. Therefore, inequalities~\eqref{eq:eberlein} and~\eqref{eq:delta} imply
\begin{equation}\label{eq:c}
\lim_{k\to\infty}\tilde{c}_{pq}^{(k)}=0.
\end{equation}
\item[(iii)] (Hari~\cite{Hari82}) For $\widetilde{A}^{(k)}=R_k^*A^{(k)}R_k$, $k\geq0$ and
\begin{equation}\label{eq:E}
E^{(k)}=A^{(k+1)}-\widetilde{A}^{(k)},
\end{equation}
we have
\begin{equation}\label{eq:hariE}
\|E^{(k)}\|_F^2\leq\frac{3}{2}n^2|\tilde{c}_{pq}^{(k)}|.
\end{equation}
\item[(iv)] (Hari~\cite{Hari82}) For $\widetilde{B}^{(k)}=R_k^*B^{(k)}R_k$, $k\geq0$ and
\begin{equation}\label{eq:F}
F^{(k)}=B^{(k+1)}-\widetilde{B}^{(k)},
\end{equation}
we have
\begin{equation}\label{eq:hariF}
\|F^{(k)}\|_F^2\leq\frac{3}{2}n^2|\tilde{c}_{pq}^{(k)}|.
\end{equation}
\item[(v)] For any $k\geq0$ we have
\begin{align}
C(\widetilde{A}^{(k)}) & = C(R_k^*A^{(k)}R_k) \nonumber \\
& = R_k^*A^{(k)}(A^{(k)})^*R_k - R_k^*(A^{(k)})^*A^{(k)}R_k \nonumber \\
& = R_k^*(A^{(k)}(A^{(k)})^*-(A^{(k)})^*A^{(k)})R_k \nonumber \\
& = R_k^*C(A^{(k)})R_k. \label{eq:CtildeA}
\end{align}
\end{itemize}

We define the \emph{off-norm} of an $n\times n$ matrix $X$ as the Frobenius norm of its off-diagonal part, that is,
$$\off^2(X)=\sum_{\substack{i,j=1 \\ i\neq j}}^n x_{ij}^2.$$
Matrix $X$ is diagonal if and only if $\off(X)=0$.

We will also use a result from~\cite{BH21} for the complex Jacobi operators. Jacobi annihilators and operators were introduced in~\cite{HZ68} and later generalized in~\cite{Hari15}. Here we give a simplified definition of the complex Jacobi annihilators and operators, the one designed to meet our needs.

For an $n\times n$ matrix $B$ we define its vectorization as a vector $b=\ve(B)\in\C^{2N}$, $N=n(n-1)/2$, containing all off-diagonal elements of $B$.
We observe a Hermitian $n\times n$ matrix $B$. Let $R$ be an $n\times n$ rotation matrix that differs from the identity matrix in its $2\times2$ submatrix $\widehat{R}$ defined by the pivot position $(p,q)$, as in~\eqref{eq:RS}, such that the rotation angle $\phi$ satisfies $|\phi|\leq\frac{\pi}{4}$. Moreover, let $\calN_{pq}:\R^{n\times n}\mapsto\R^{n\times n}$ be an operator that sets to zero the elements on positions $(p,q)$ and $(q,p)$ in a given matrix. A \emph{complex Jacobi annihilator} $\calR_{pq}(R)\in\C^{2N\times2N}$ is defined by the rule
\begin{equation}\label{Jannihilator}
\calR_{pq}(R)\ve(B):=\ve(\calN_{pq}(R^*BR)).
\end{equation}

For a pivot ordering $\calo=(p_0,q_0),(p_1,q_1),\ldots,(p_{N-1},q_{N-1})\in\calO(\calP_n)$, a \emph{complex Jacobi operator} determined by the ordering $\calo$ is defined as a product of $N$ Jacobi annihilators
\begin{equation*}
\calJ_{\calo}:=\calR_{p_{N-1},q_{N-1}}(R_{N-1})\cdots\calR_{p_1,q_1}(R_1)\calR_{p_0,q_0}(R_0).
\end{equation*}

The definitions presented above are the special cases of the prevailing definitions from~\cite{BH21}. It will be useful to know that the spectral norm of a Jacobi annihilator is equal to one, except for the case of the $2\times2$ annihilator which is a zero-matrix. This follows from the structure of the annihilators (see, e.g.~\cite{BH21}).

\begin{Proposition}\label{prop_sgJ}
Let $\calo\in\calC^{(n)}_{sg}$. Suppose that
$\calo\stackrel{p}{\sim}\calo'\stackrel{w}{\sim}\calo''$ or $\calo\stackrel{w}{\sim}\calo'\stackrel{p}{\sim}\calo''$, $\calo''\in\calC^{(n)}_{sp}$, and that the weak equivalence relation is in the canonical form containing exactly $d$ shift equivalences.
Then, for any $d+1$ Jacobi operators $\calJ_{\calo,1},\calJ_{\calo,2},\ldots,\calJ_{\calo,d+1}$, there is a constant $\gamma_n$ depending only on $n$ such that
$$\|\calJ_{\calo,1}\calJ_{\calo,2}\cdots\calJ_{\calo,d+1}\|_2\leq\gamma_n, \quad 0\leq\gamma_n<1.$$
\end{Proposition}

\begin{proof}
This is a special case of Theorem 3.6 from~\cite{BH21}.
\end{proof}

Further on, we prove the following two auxiliary propositions.

\begin{Proposition}\label{lemma:sequence}
Let $(x_k,k\geq0)$ be a sequence of nonnegative real numbers such that
\begin{equation}\label{lemma:xk}
x_{k+1}=\gamma x_k+c_k, \quad 0\leq\gamma<1.
\end{equation}
If $\lim_{k\to\infty}c_k=0$, then
$$\lim_{k\to\infty}x_k=0.$$
\end{Proposition}

\begin{proof}
First, we show that the sequence~\eqref{lemma:xk} is bounded from above. Take
$$C=\max\{x_0,\sup_k c_k\}.$$
We prove the boundedness by mathematical induction. For $k=0$,
$$x_0\leq C\leq\frac{C}{1-\gamma}=:M, \quad \text{for } 0\leq\gamma<1.$$
Assume that $x_k\leq M$ for some given $k$. Then, for $k+1$,
$$x_{k+1}=\gamma x_k+c_k \leq\gamma M+C=\gamma M+(1-\gamma)M=M.$$
Therefore, $x_k\leq M$ for any $k\geq0$.

For the limit superior we take $\limsup_{k\to\infty}x_k=L\in\R$. Then,
$$L=\limsup_{k\to\infty}x_{k+1} \leq\gamma\limsup_{k\to\infty}x_k+\limsup_{k\to\infty}c_k=\gamma L.$$
Since $0\leq\gamma<1$, the upper inequality can hold only with $L=0$. Since $(x_k)_k$ is the sequence of nonnegative real numbers, $\liminf_{k\to\infty}x_k\geq0$. This implies that
$$\limsup_{k\to\infty}x_k=\liminf_{k\to\infty}x_k=0$$
and $\lim_{k\to\infty}x_k=0.$
\end{proof}

\begin{Proposition}\label{prop:sg}
Let $H\neq 0$ be a Hermitian matrix. Let $(H^{(k)}, k\geq 0)$ be a sequence generated by applying the following iterative process on $H$,
\begin{equation}\label{prop:Hprocess}
H^{(k+1)}=R_k^*H^{(k)}R_k + M^{(k)}, \quad H^{(0)}=H,\quad k\geq 0,
\end{equation}
where  $R_k$ are complex plane rotations acting in the $(p(k),q(k))$ plane, $p(k)<q(k)$, with the rotation angles $\left|\varphi_k\right|\leq \frac{\pi}{4}$, $k\geq 0$. Suppose that the pivot strategy is defined by an ordering $\calo\in\calC^{(n)}_{sg}$ and that
\begin{equation}\label{prop:offM}
\lim_{k\to\infty} \off(M^{(k)})=0.
\end{equation}
Then the limit
\begin{equation}\label{prop:limh}
\lim_{k\to\infty} \left| h_{p(k) q(k)}^{(k+1)} \right|=0,
\end{equation}
implies
$$\lim_{k\to\infty}\off(H^{(k)})=0.$$
\end{Proposition}

\begin{proof}
The proof uses the idea of the proof of Theorem 3.8 from~\cite{BH21}.

To simplify the notation, let $(p,q)=(p(k),q(k))$ denote the pivot pair at step $k$.
Transformation $R_k^*H^{(k)}R_k$ does not annihilate the elements on positions $(p,q)$ and $(q,p)$ of $H^{(k)}$, but we can write it as
\begin{equation}\label{prop:N}
R_k^*H^{(k)}R_k=\calN_{pq}(R_k^*H^{(k)}R_k)+(R_k^*H^{(k)}R_k)_{pq}(e_{p}e_{q}^*)+(R_k^*H^{(k)}R_k)_{qp}(e_{q}e_{p}^*),
\end{equation}
where $e_r$ is the $r$th column vector of the identity matrix $I_n$ and $\calN_{pq}$ is as in~\eqref{Jannihilator}.
By using the $\ve$ operator on equation~\eqref{prop:Hprocess} and the definition of a Jacobi annihilator~\eqref{Jannihilator}, from the relation~\eqref{prop:N} we get
\begin{equation}\label{prop:chi}
\chi^{(k+1)}=\calR_{p_k q_k}(R_k) \chi^{(k)}+m^{(k)}, \quad k\geq 0,
\end{equation}
where $\chi^{(k)}=\ve(H^{(k)})$, and
\begin{equation}\label{prop:m}
m^{(k)}=\ve(M^{(k)})+(R_k^*H^{(k)}R_k)_{pq}e_{\tau(p,q)}+(R_k^*H^{(k)}R_k)_{qp}e_{\tau(q,p)}.
\end{equation}
Here, $\tau(p,q)$ stands for the position of the matrix element $x_{pq}$ in the vectorization $\ve(X)$ and $e_{\tau(p,q)}$ is the column vector of the identity matrix $I_{2N}$ with one on position $\tau(p,q)$.
It is easy to see that the relation~\eqref{prop:m} together with the assumptions~\eqref{prop:offM} and~\eqref{prop:limh} imply that
\begin{equation}\label{prop:lim_mk}
\lim_{k\to \infty}m^{(k)}=0.
\end{equation}

We denote the matrix obtained from $H$ after $t$ cycles of the process~\eqref{prop:Hprocess} by $H^{(tN)}$.
Vector $\chi^{(tN)}=\ve(H^{(tN)})$ can be written as
$$\chi^{(tN)}=\calJ_{\calo}^{[tN]} \chi^{((t-1)N)}+m^{[tN]}, \quad t\geq 1,$$
The Jacobi operator $\calJ_{\calo}^{[tN]}$ that appears in the upper equation
is determined by the ordering $\calo=(p_0,q_0),(p_1,q_1),\ldots,(p_{N-1},q_{N-1})\in\calO(\calP_n)$ and by the Jacobi annihilators,
$$\calJ_{\calo}^{[tN]}=\calR_{p_{N-1},q_{N-1}}(R_{tN-1})\cdots\calR_{p_1,q_1}(R_{(t-1)N+1})\calR_{p_0,q_0}(R_{(t-1)N}),$$
while
\begin{align}
m^{[tN]} & = \calR_{p_{N-1},q_{N-1}}(R_{tN-1})\cdots\calR_{p_1,q_1}(R_{(t-1)N+1})\calR_{p_0,q_0}(R_{(t-1)N})m^{((t-1)N)} \nonumber \\
& \qquad +\cdots+\calR_{p_{N-1},q_{N-1}}(R_{tN-1})m^{(tN-2)}+m^{(tN-1)}. \label{prop:mtN}
\end{align}
From the fact that the spectral norm of any Jacobi annihilator is equal to one (or zero if it is a $2\times2$ annihilator), the relation~\eqref{prop:mtN} indicates that
$$\|m^{[tN]}\|_2 \leq \| m^{((t-1)N)}\|_2+\cdots+\|m^{(tN-2)}\|_2+\| m^{(tN-1)}\|_2, \quad t\geq1.$$
Thus, from the limit~\eqref{prop:lim_mk} we get
\begin{equation} \label{prop:lim_mtN}
\lim_{t\to\infty}m^{[tN]}=0.
\end{equation}

Since $\calo\in \calC_{sg}^{(n)}$, i.e.\@ the pivot strategy is generalized serial, suppose that
$\calo\stackrel{p}{\sim}\calo'\stackrel{w}{\sim}\calo''$ or $\calo\stackrel{w}{\sim}\calo'\stackrel{p}{\sim}\calo''$, $\calo''\in\calC^{(n)}_{sp}$, and that the weak equivalence relation is in the canonical form containing exactly $d$ shift equivalences. For $d+1$ consecutive cycles we get
\begin{equation}\label{prop:chi_t+dN}
\chi^{((t+d)N)} = \calJ_{\calo}^{[(t+d)N]}\cdots\calJ_{\calo}^{[(t+1)N]}\calJ_{\calo}^{[tN]}\chi^{((t-1)N)}+m_{[d+1]}^{[tN]}, \quad t\geq 1,
\end{equation}
where
\begin{equation*}
m_{[d+1]}^{[tN]} = \calJ_{\calo}^{[(t+d)N]}\cdots\calJ_{\calo}^{[(t+1)N]}m^{[tN]} +\cdots+\calJ_{\calo}^{[(t+d)N]}m^{[(t+d-1)N]} +m^{[(t+d)N]}.
\end{equation*}
Similarly as before, the property of the spectral norm of the Jacobi operator implies
$$\|m_{[d+1]}^{[tN]}\|_2 \leq \| m^{[tN]}\|_2+\| m^{[(t+1)N]}\|_2+\cdots+\| m^{[(t+d)N]}\|_2,$$
and using the limit~\eqref{prop:lim_mtN} we get
\begin{equation*}
\lim_{t\to\infty}m_{[d+1]}^{[tN]}=0.
\end{equation*}
To the Jacobi operators from~\eqref{prop:chi_t+dN} we can apply the Proposition~\ref{prop_sgJ}. We get
\begin{equation}~\label{prop:Jgamma}
\|\calJ_{\calo}^{[(t+d)N]}\cdots\calJ_{\calo}^{[(t+1)N]}\calJ_{\calo}^{[tN]}\|_2\leq\gamma_n, \quad 0\leq\gamma_n<1.
\end{equation}
Looking at the spectral norm of~\eqref{prop:chi_t+dN} and using the bound~\eqref{prop:Jgamma} we obtain
\begin{align*}
\|\chi^{[(t+d)N)}\|_2 & \leq \|\calJ_{\calo}^{[(t+d)N]}\cdots\calJ_{\calo}^{[(t+1)N]}\calJ_{\calo}^{[tN]}\|_2\|\chi^{((t-1)N)}\|_2 + \|m_{[d+1]}^{[tN]}\|_2 \\
& \leq \gamma_n\|\chi^{[(t-1)N]}\|_2 + \|m_{[d+1]}^{[tN]}\|_2.
\end{align*}
Considering that $0\leq \gamma_n<1$ and $\|m_{[d+1]}^{[tN]}\|_2\to0$, as $t\to\infty$, we employ the Proposition~\ref{lemma:sequence} which yields $\lim_{t\to\infty}\chi^{(tN)}=0$. Therefore, iterations obtained after each cycle converge to zero.

Additionally, for iterations $0<k<N$ within one cycle, from the relation~\eqref{prop:chi} we have
\begin{align*}
\chi^{((t-1)N+k)} & =\calR_{p_{k-1},q_{k-1}}(R_{(t-1)N+k-1})\cdots\calR_{p_1,q_1}(R_{(t-1)N+1})\calR_{p_0,q_0}(R_{(t-1)N})\chi^{((t-1)N)} \\
& \qquad + \calR_{p_{k-1},q_{k-1}}(R_{(t-1)N+k-1})\cdots\calR_{p_1,q_1}(R_{(t-1)N+1})\calR_{p_0,q_0}(R_{(t-1)N})m^{((t-1)N)} \\
& \qquad +\cdots+\calR_{p_{k-1},q_{k-1}}(R_{(t-1)N+k-1})m^{((t-1)N+k-2)}+m^{((t-1)N+k-1)}.
\end{align*}
In the same manner as before we get the inequality
\begin{align*}
\|\chi^{((t-1)N+k)}\|_2 & \leq \|\chi^{((t-1)N)}\|_2+\| m^{((t-1)N)}\|_2+\cdots+\|m^{((t-1)N+k-2)}\|_2+\|m^{((t-1)N+k-1)}\|_2 \\
& \leq \|\chi^{((t-1)N)}\|_2 + k\max_{0\leq r\leq k-1} \|m^{((t-1)N+r)}\|_2.
\end{align*}
Thus, $\lim_{t\to\infty}\|\chi^{((t-1)N+k)}\|_2=0$,
and it follows $$\lim_{k\to\infty} \|\chi^{(k)}\|_2=0.$$
Finally, because $\off(H^{(k)})=\|\chi^{(k)}\|_2$, $k\geq0$, we have $\lim_{k\to\infty}\off(H^{(k)})=0.$
\end{proof}

Now we can prove the convergence theorem for Eberlein method under the serial orderings with permutations, $\calo\in\calC_{sg}^{(n)}$. Matrix $B^{(k)}$ is defined as in the equation~\eqref{eq:B} and $C(B^{(k)})$ is as in~\eqref{eq:C}.

\begin{Theorem}\label{theorem:sp}
Let $A\in\C^{n\times n}$ and let $(A^{(k)},k\geq0)$  be a sequence generated by the Eberlein method under a generalized serial pivot strategy defined by an ordering $\calo\in\calC^{(n)}_{sg}$. Then
\begin{itemize}
\item[(i)] The sequence of the off-norms $(\off(B^{(k)}),k\geq0)$ tends to zero,
$$\lim_{k\to\infty}\off(B^{(k)})=0.$$
\item[(ii)] The sequence $(A^{(k)},k\geq0)$ tends to a normal matrix, that is,
$$\lim_{k\to\infty}C(A^{(k)})=0.$$
\item[(iii)] The sequence of matrices $(B^{(k)},k\geq0)$ tends to a fixed diagonal matrix,
$$\lim_{k\to\infty}B^{(k)}=\diag(\mu_1,\mu_2,\ldots,\mu_n),$$
where $\mu_i$, $1\leq i\leq n$, are real parts of the eigenvalues of $A$.
\item[(iv)] If $\mu_i\neq\mu_j$, then $\lim_{k\to\infty}a_{ij}^{(k)}=0$ and $\lim_{k\to\infty}a_{ji}^{(k)}=0$.
\end{itemize}
\end{Theorem}

\begin{proof}
\begin{itemize}
\item[(i)]
For $F^{(k)}$ defined as in~\eqref{eq:F} we have
\begin{equation}\label{tm:B}
B^{(k+1)}=R_k^*B^{(k)}R_k+F^{(k)}, \quad k\geq0.
\end{equation}
On the pivot position $(p,q)$ in the step $k$ we have
$$b_{pq}^{(k+1)}=\tilde{b}_{pq}^{(k)}+f_{pq}^{(k)},$$
where $F^{(k)}=(f_{ij}^{(k)})$.

Relations~\eqref{eq:hariF} and~\eqref{eq:c} imply $\lim_{k\to\infty}F^{(k)}=0$ and $\lim_{k\to\infty}f_{pq}^{(k)}=0$.
Furthermore, the rotation $R_k$ is chosen to annihilate $\tilde{b}_{pq}^{(k)}$. Therefore, $\lim_{k\to\infty}b_{pq}^{(k+1)}=0$.
Matrix $B^{(0)}=B$ is Hermitian by the definition and the iterative process~\eqref{tm:B} satisfies the assumptions of the Proposition~\ref{prop:sg}. Hence,
\begin{equation*}
\lim_{k\to\infty}\off(B^{(k)})=0.
\end{equation*}

\item[(ii)]
For $E^{(k)}$ defined as in~\eqref{eq:E} we have
$$C(A^{(k+1)})=C(\tilde{A}^{(k)}+E^{(k)}).$$
Then,
\begin{align}
C(A^{(k+1)}) & = (\tilde{A}^{(k)}+E^{(k)})(\tilde{A}^{(k)}+E^{(k)})^*-(\tilde{A}^{(k)}+E^{(k)})^*(\tilde{A}^{(k)}+E^{(k)}) \nonumber \\
& = \tilde{A}^{(k)}(\tilde{A}^{(k)})^*+E^{(k)}(\tilde{A}^{(k)})^*+(\tilde{A}^{(k)}+E^{(k)})(E^{(k)})^* \nonumber \\
& \qquad - (\tilde{A}^{(k)})^*\tilde{A}^{(k)}-(E^{(k)})^*\tilde{A}^{(k)}-(\tilde{A}^{(k)}+E^{(k)})^*E^{(k)} \nonumber \\
& = C(\tilde{A}^{(k)})+A^{(k+1)}(E^{(k)})^*-(A^{(k+1)})^*E^{(k)}+E^{(k)}(\tilde{A}^{(k)})^*-(E^{(k)})^*\tilde{A}^{(k)} \nonumber \\
& = C(\tilde{A}^{(k)})+W^{(k)}, \label{tm:CA}
\end{align}
where
$$W^{(k)}=A^{(k+1)}(E^{(k)})^*-(A^{(k+1)})^*E^{(k)}+E^{(k)}(\tilde{A}^{(k)})^*-(E^{(k)})^*\tilde{A}^{(k)}.$$
Moreover, applying the relation~\eqref{eq:CtildeA}, we can write~\eqref{tm:CA} as
\begin{equation*}
C(A^{(k+1)})=R_k^*C(A^{(k)})R_k+W^{(k)}.
\end{equation*}
Using the properties of the norm and the inequality~\eqref{eq:eberlein} we get
\begin{align*}
\|W^{(k)}\|_F & \leq \|A^{(k+1)}(E^{(k)})^*\|_F+\|(A^{(k+1)})^*E^{(k)}\|_F+\|E^{(k)}(\tilde{A}^{(k)})^*\|_F+\|(E^{(k)})^*\tilde{A}^{(k)}\|_F \\
& = 2\|E^{(k)}\|_F(\|A^{(k+1)}\|_F+\|\tilde{A}^{(k)}\|_F) \\
& \leq 4\|E^{(k)}\|_F\|\tilde{A}^{(k)}\|_F,
\end{align*}
and
$$\|W^{(k)}\|_F^2\leq 16\|E^{(k)}\|_F^2\|\tilde{A}^{(k)}\|_F^2.$$
It follows from the relations~\eqref{eq:E} and~\eqref{eq:hariE} that
$$\|W^{(k)}\|_F^2\leq 16\|E^{(k)}\|_F^2\|A\|_F^2\leq24n^2|\tilde{c}_{pq}^{(k)}|\|A\|_F^2.$$
Thus, relation~\eqref{eq:c} implies
\begin{equation}\label{tm:W}
\lim_{k\to\infty}\|W^{(k)}\|_F=0.
\end{equation}

We consider the off-diagonal and the diagonal part of $C(A^{(k)})$ separately.
Similarly as for matrices $B^{(k)}$, on the pivot position $(p,q)$ in the step $k$ we have
$$c_{pq}^{(k+1)}=\tilde{c}_{pq}^{(k)}+w_{pq}^{(k)},$$
where $W^{(k)}=(w_{ij}^{(k)})$. Relations~\eqref{eq:c} and~\eqref{tm:W} imply $\lim_{k\to\infty}c_{pq}^{(k+1)}=0$. It is easy to check that matrices $C(A^{(k)})$, $k\geq0$, are Hermitian and we can use the Proposition~\ref{prop:sg} again. We get
\begin{equation}\label{tm:offC}
\lim_{k\to\infty}\off(C(A^{(k)}))=0.
\end{equation}

It remains to show that
$$\lim_{k\to\infty}c_{ii}^{(k)}=0.$$
Set $A^{(k)}=B^{(k)}+Z^{(k)}$, where $(B^{(k)})$ is Hermitian, as in~\eqref{eq:B}, and $Z^{(k)}$ is skew-Hermitian. Then,
\begin{equation}\label{tm:CBG}
C(A^{(k)})=2(Z^{(k)}B^{(k)}-B^{(k)}Z^{(k)}).
\end{equation}
The diagonal element of $C(A^{(k)})$ is given by
$$c_{ii}^{(k)}=2\sum_{j=1}^n \left(z_{ij}^{(k)} b_{ji}^{(k)}-b_{ij}^{(k)}z_{ji}^{(k)}\right).$$
It is proven in part $(i)$ that $\lim_{k\to \infty} \off(B^{(k)})=0$, that is,
$$\lim_{k\to \infty} b_{ij}^{(k)}=0, \quad \text{for } i\neq j.$$
Thus,
\begin{equation}\label{tm:cii}
\lim_{k\to\infty}c_{ii}^{(k)}=2\left(z_{ii}^{(k)}b_{ii}^{(k)}-b_{ii}^{(k)}z_{ii}^{(k)}\right)=0.
\end{equation}
Relations~\eqref{tm:offC} and~\eqref{tm:cii} imply the assertion $(ii)$ of the theorem.

\item[(iii)]
In part $(i)$ of the proof we showed that matrices $B^{(k)}$ tend to a diagonal matrix. The fact that the diagonal elements of the matrix $\lim_{k\to\infty}B^{(k)}$ correspond to the real parts of the eigenvalues of $A$ is then proved as in~\cite{PupHari99}, using the assertion $(ii)$ of this theorem.

\item[(iv)]
Using the relation~\eqref{tm:CBG} and parts $(i)$ and $(ii)$ of the theorem it follows that
\begin{align*}
0 & = \lim_{k\to\infty}c_{ij}^{(k)} =2\lim_{k\to \infty}\sum_{r=1}^n \left(z_{ir}^{(k)}b_{rj}^{(k)}-b_{ir}^{(k)}z_{rj}^{(k)}\right) \\
& = 2\lim_{k\to\infty}\left(z_{ij}^{(k)}b_{jj}^{(k)}-b_{ii}^{(k)}z_{ij}^{(k)}\right) =2(\mu_j-\mu_i) \lim_{k\to\infty}z_{ij}^{(k)}, \quad 1\leq i,j\leq n.
\end{align*}
If $\mu_i\neq\mu_j$, then $\lim_{k\to\infty}z_{ij}^{(k)}=0$. Finally, since $\lim_{k\to\infty}b_{ij}^{(k)}=0$ for $i\neq j$, we get
\begin{equation*}
a_{ij}^{(k)}=b_{ij}^{(k)}+z_{ij}^{(k)}\to0 \quad \text{and} \quad a_{ji}^{(k)}=(b_{ij}^{(k)})^*-(z_{ij}^{(k)})^*\to 0, \quad \text{as } k\to\infty.
\end{equation*}
\end{itemize}
\end{proof}

Therefore, starting with an $n\times n$ matrix $A$, the Eberlein method under pivot strategies defined by any generalized serial pivot strategy converges to some matrix $\Lambda$. Assuming that the diagonal elements of $\Lambda$ are arranged such that their real parts appear in decreasing order, based on Theorem~\ref{theorem:sp}, we get to the following conclusion. Matrix $\Lambda$ is a block diagonal matrix with block sizes corresponding to the multiplicities of the real parts of the eigenvalues of $A$. In order to find all eigenvalues of $A$, it remains to find the eigenvalues of the diagonal blocks. To that end one can, e.g., apply the nonsymmetric Jacobi algorithm for the computation of the Schur form discussed in~\cite{Mehl08}.

\section{Numerical results}\label{sec:num}

Numerical tests of the Algorithm~\ref{agm:eberlein} under the generalized pivot strategies with permutations are presented in this section. All experiments are done in Matlab R2021a.

\begin{figure}[!h]
\begin{subfigure}{0.45\textwidth}
    \centering
    \includegraphics[width=\textwidth]{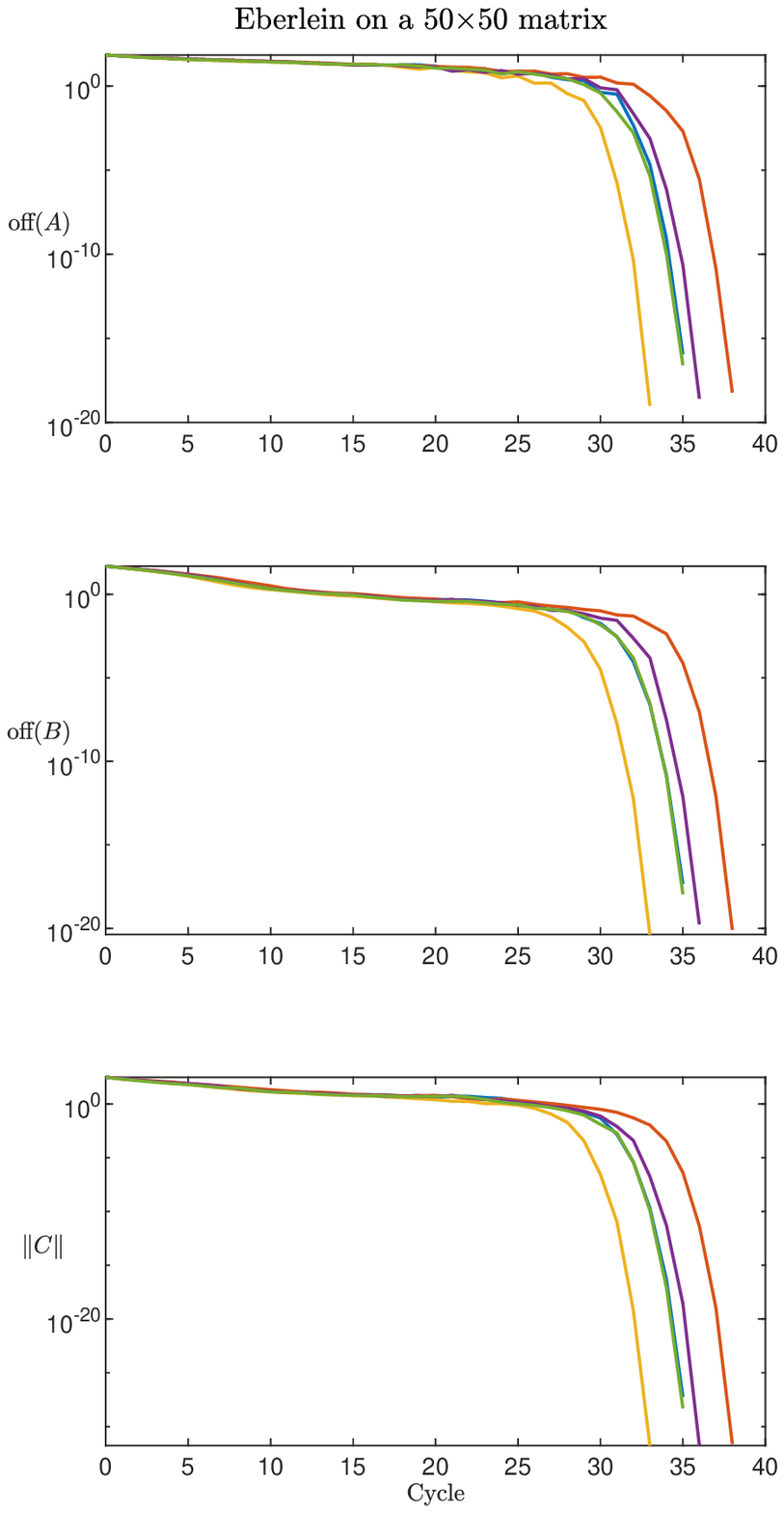}
    \caption{Complex algorithm, random $A\in\C^{50\times 50}$}
\end{subfigure}
\hfill
\begin{subfigure}{0.45\textwidth}
        \includegraphics[width=\textwidth]{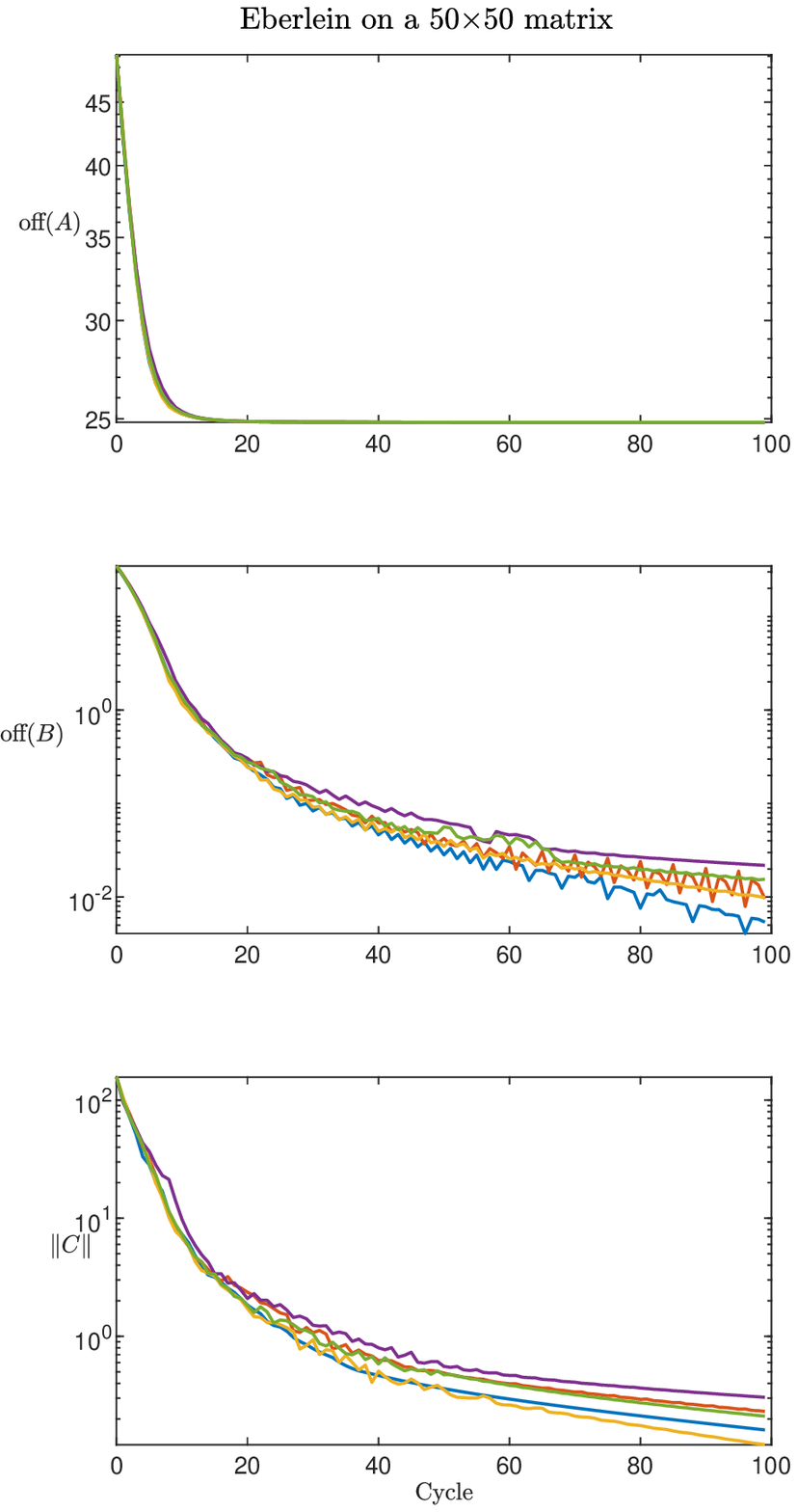}
    \caption{Real algorithm, random $A\in\R^{50\times 50}$}
\end{subfigure}
\caption{Change in $\off(A^{(k)})$, $\off(B^{(k)})$ and $\|C(A^{(k)})\|_F$.}
\label{fig:n50}
\end{figure}

To depict the performance of the Eberlein algorithm, we observe three quantities; $\off(A^{(k)})$, $\off(B^{(k)})$, and $\|C(A^{(k)})\|_F$. The results are presented in logarithmic scale. The algorithm is terminated when the change in the off-norm of $B^{(k)}$ becomes small enough, $10^{-8}$. According to Theorem~\ref{theorem:sp}, both $\off(B^{(k)})$ and $\|C(A^{(k)})\|_F$ should converge to zero.
In Figure~\ref{fig:n50} the results of the Eberlein algorithm on a non-structured random complex matrix are shown, as well as the results of the real Eberlein algorithm on a non-structured random real matrix.
Each line represents a different pivot strategy $I_\calo$, $\calo\in\calC^{(n)}_{sg}$, chosen randomly at the beginning of the algorithm.

\begin{figure}[!h]
    \centering
    \includegraphics[width=\textwidth]{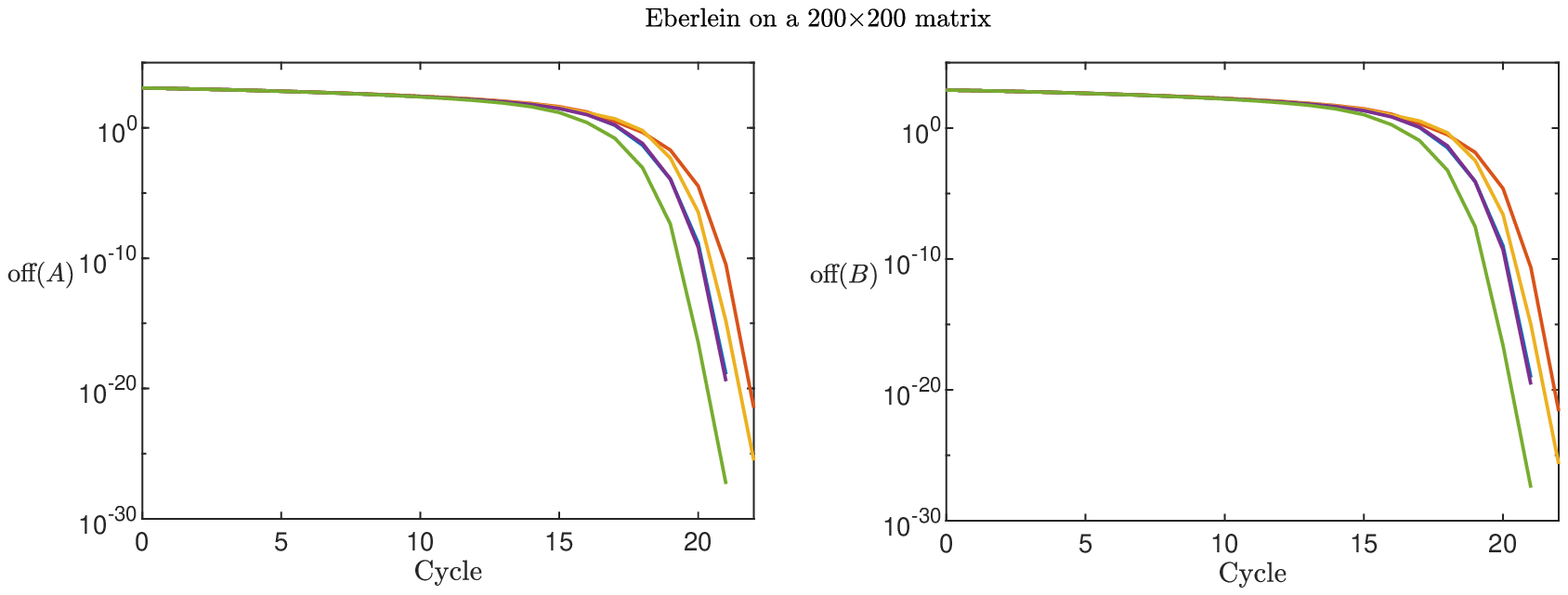}
    \caption{Progress of the off-norms of $A^{(k)}$ and $B^{(k)}$, in logarithmic scale, for a unitarily diagonalizable matrix.}
    \label{fig:n200cplx}
\end{figure}

The algorithm is significantly faster if it is applied on a normal matrix.
We construct a unitarily diagonalizable matrix $A=A^{(0)}$ such that we multiply some chosen complex diagonal matrix from the left and right hand side by a random unitary matrix.
In Figure~\ref{fig:n200cplx}, we see the results of the Eberlein method applied on a diagonalizable complex matrix. Here we do not show $\|C(A^{(k)})\|$ because $A^{(0)}$ is normal, that is $C(A^{(0)})=0$, and it stays normal during the process.

\begin{figure}[!h]
\begin{subfigure}[t]{\textwidth}
    \includegraphics[width=\textwidth]{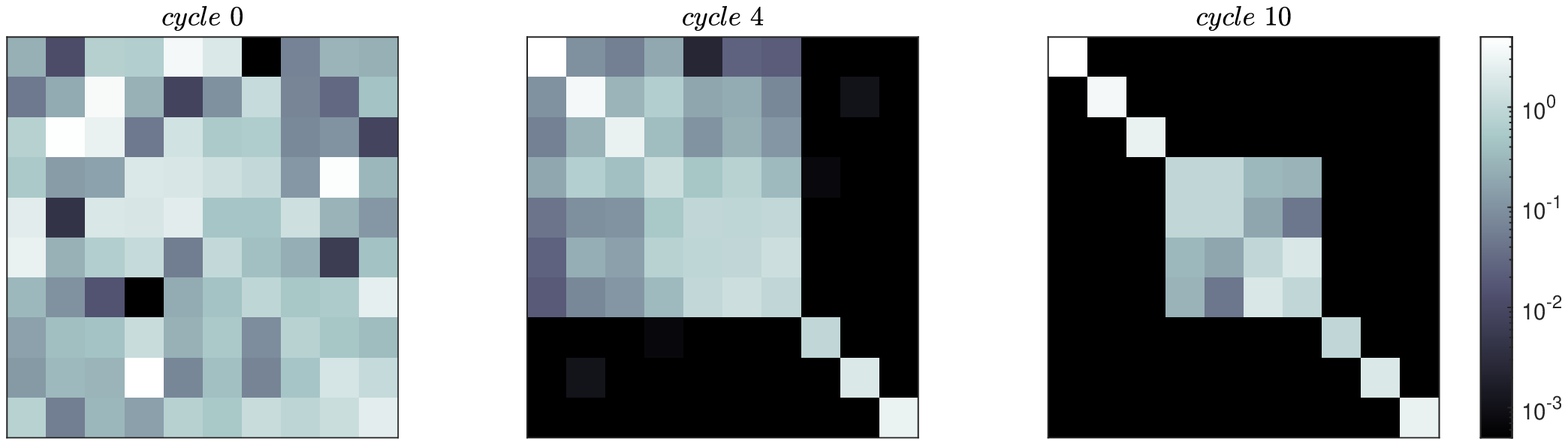}
    \caption{Two complex conjugate pairs of eigenvalues with the same real part that formed a $4\times 4$ diagonal block. }
    \label{subfig1}
\end{subfigure}\hspace{\fill} 
\bigskip 
\begin{subfigure}[t]{\textwidth}
    \includegraphics[width=\linewidth]{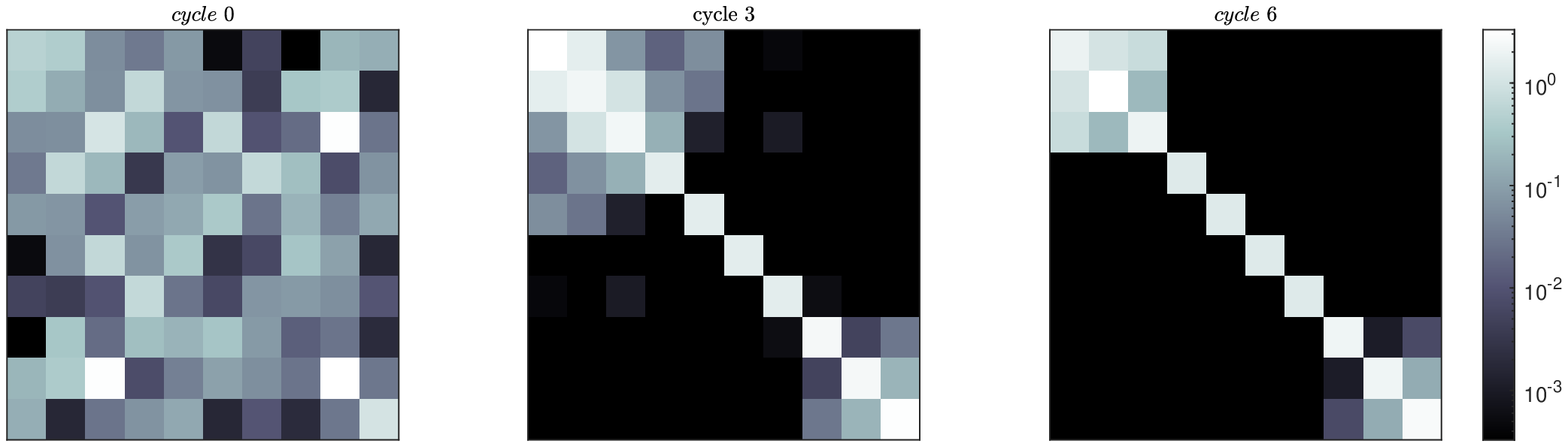}
    \caption{Different eigenvalues with the same real part formed $3\times 3$ diagonal blocks, while the rest of the diagonal carries the quadruple eigenvalue with unique real part.}
    \label{subfig2}
\end{subfigure}
\caption{Block diagonal structure.}
\label{fig:heat}
\end{figure}

In order to show the block diagonal structure of $A^{(k)}$ discussed at the end of the previous section, we applied the Eberlein method on two matrices from $\C^{10\times 10}$. To generate the starting matrix $A$, first we set upper-triangular matrix $T$ to have the diagonal elements listed below. Then we multiply $T$ by a random unitary matrix $Q$, $A=Q^*TQ$. In our implementation of the algorithm, we introduce an additional condition so that the real values of the diagonal elements appear in the decreasing order. That is achieved by, if necessary, translating the angle $\alpha_k$ by $\pi/2$ in the $k$th step of the process. The evolution of the matrix structure of the iterates is shown in Figure~\ref{fig:heat}. Specifically, the figure shows the logarithm of the absolute values of the elements of $A^{(k)}$. Lighter squares represent elements larger in absolute value. According to the Theorem~\ref{theorem:sp}, the algorithm should converge to a block diagonal matrix in both cases described below.

In Figure~\ref{subfig1}, matrix has distinctive eigenvalues with the spectrum
$$\{5,4,3, 1\pm 2i, 1\pm i, -1,-2,-3 \}.$$
Thus, we deal with two complex conjugate pairs of eigenvalues with the same real part.
On the other hand, in Figure~\ref{subfig2}, the matrix spectrum is
$$\{-2+i,-2+i,-2+2i, 1-i, 1-i, 1-i, 1-i, 2-i,2+3i,2+i \}.$$
The key is to observe that all eigenvalues with the real part equal to 1 have the same imaginary parts, i.e.\@ $1-i$ has algebraic multiplicity 4. On the contrary, other eigenvalues with equal real parts have at least two distinctive imaginary parts.

For both matrices, after a few cycles we can faintly see the diagonal blocks. After a few more cycles the block diagonal structure is clear. For the first matrix, the obtained $4\times 4$ block has eigenvalues that are (approximately) $1\pm i$ and $1\pm 2i$. The rest of the diagonal carries the real eigenvalues of the original matrix.
On the other hand, for the second matrix we see two blocks that correspond to eigenvalues with real parts equal to 2 and -2. The rest of the diagonal corresponds to the eigenvalue $1-i$ and it does not form a block despite the quadruple multiplicity. The difference is that there are no other eigenvalues with the same real part, but different imaginary part.


\begin{thebibliography}{10}

\bibitem{BePhD}
E. Begović:
\emph{Konvergencija blok Jacobijevih metoda}.
Ph.D. thesis, University of Zagreb, Faculty of Science, Zagreb (2014)
%
\bibitem{BGFass97}
A. Bunse-Gerstner, H. Faßbender:
\emph{A Jacobi-like method for solving algebraic Riccati equations on parallel computers}.
IEEE Trans. Automat. Control, 42 (1997) 1071--1084.
%
\bibitem{DeVe92}
J. Demmel, K. Veseli\'c:
\emph{Jacobi's method is more accurate than QR}.
SIAM J. Matrix Anal. Appl. 13(4) (1992) 1204--1245.
%
\bibitem{DV08-1}
Z. Drma\v c, K. Veseli\'c:
\emph{New fast and accurate Jacobi SVD algorithm I}.
SIAM J. Matrix Anal. Appl. 29(4) (2008) 1322--1342.
%
\bibitem{DV08-2}
Z. Drma\v c, K. Veseli\'c:
\emph{New fast and accurate Jacobi SVD algorithm II}.
SIAM J. Matrix Anal. Appl. 29(4) (2008) 1343--1362.
%
\bibitem{Eber62}
P. J. Eberlein:
\emph{A Jacobi-like method for the automatic computation of eigenvalues and eigenvectors of an arbitrary matrix}.
SIAM J. 10(1) (1962) 74--88.
%
\bibitem{FMM01}
H. Fa{\ss}bender, D. S. Mackey, N. Mackey:
\emph{Hamiltonian and Jacobi come full circle: Jacobi algorithms for structured Hamiltonian eigenproblems}.
Linear Algebra Appl., 332--334 (2001) 37--80.
%
\bibitem{Hari82}
V. Hari:
\emph{On the global convergence of the Eberlein method for real matrices}.
Numer. Math. 39 (1982) 361--369.
%
\bibitem{Hari15}
V. Hari:
\emph{Convergence to diagonal form of block Jacobi-type methods}.
Numer. Math. 129(3) (2015) 449--481.
%
\bibitem{BH17}
V. Hari, E. Begovi\'{c}~Kova\v{c}:
\emph{Convergence of the cyclic and quasi-cyclic block Jacobi methods}.
Electron. Trans. Numer. Anal. 46 (2017) 107-147.
%
\bibitem{BH21}
V. Hari, E. Begović Kovač:
\emph{On the convergence of complex Jacobi methods}.
Linear multilinear algebra. 69(3) (2021) 489--514.
%
\bibitem{HSS14}
V. Hari, S. Singer, S. Singer:
\emph{Full block $J$--Jacobi method for Hermitian matrices}.
Linear Algebra Appl. 444 (2014) 1--27.
%
\bibitem{HZ68}
P. Henrici, K. Zimmermann:
\emph{An estimate for the norms of certain cyclic Jacobi operators}.
Linear Algebra Appl. 1(4) (1968) 489–501.
%
\bibitem{MMMM09}
D. S. Mackey, N. Mackey, C. Mehl, V. Mehrmann:
\emph{Numerical methods for palindromic eigenvalue problems: Computing the anti-triangular Schur form}.
Numer. Linear Algebra Appl. 16(1) (2009) 63--86.
%
\bibitem{MMT03}
D. S. Mackey, N. Mackey, F. Tisseur:
\emph{Structured tools for structured matrices}.
Electron. J. Linear Algebra 10 (2003) 106--145.
%
\bibitem{Masc95}
W. F. Mascarenhas:
\emph{On the convergence of the Jacobi method for arbitrary orderings}.
SIAM J. Matrix Anal. Appl. 16(4) (1995) 1197--1209.
%
\bibitem{Matejas09}
J. Mateja\v s:
\emph{Accuracy of the Jacobi method on scaled diagonally dominant symmetric matrices}.
SIAM J. Matrix Anal. Appl. 31(1) (2009) 133--153.
%
\bibitem{Mehl04}
C. Mehl:
\emph{Jacobi-like algorithms for the indefinite generalized Hermitian eigenvalue problem}.
SIAM J. Matrix Anal. Appl., 25 (2004) 964--985.
%
\bibitem{Mehl08}
C. Mehl: \emph{On asymptotic convergence of nonsymmetric Jacobi algorithms}.
SIAM J. Matrix Anal. Appl. 30(1) (2008) 291--311.
%
\bibitem{PupHari99}
D. Pupovci, V. Hari:
\emph{On the convergence of parallelized Eberlein methods}.
Rad. Mat. 8 (1992) 249--267.
%
\bibitem{Ves76}
K. Veselić:
\emph{A convergent Jacobi method for solving the eigenproblem of arbitrary real matrices}.
Numer. Math. 25 (1976) 179--184.
\end{thebibliography}
\end{document}